\documentclass[11pt, a4paper,twoside]{amsart}

  \usepackage[margin=2.5cm]{geometry}
    \usepackage{multicol}
  \usepackage[T1]{fontenc}
  \usepackage[utf8]{inputenc}
  \usepackage[british]{babel}
    \usepackage{microtype}

  \usepackage{newtxtext,newtxmath}
  \usepackage{helvet}
  \usepackage{xcolor}
    \definecolor{Navy}{HTML}{12355B}
    \definecolor{DarkGray}{HTML}{222222}
    \definecolor{LightGray}{HTML}{666666}

\usepackage{url}
\usepackage[normalem]{ulem}
  \usepackage{threeparttable}

\color{DarkGray}  

\usepackage[explicit]{titlesec}

\titleformat{\section}
{\Large\sffamily\bfseries\color{Navy}}
{\thesection}
{0.8em}
{#1}
[\vspace{0.3em}\titlerule]

\titleformat{\subsection}
{\large\sffamily\bfseries\color{Navy}}
{\thesubsection}
{0.7em}
{#1 \vspace{-.7em}}

\titlespacing*{\section}
{0pt}{3ex}{1.2ex}

\makeatletter  
  \titleformat{\subsubsection}[runin]
  {\bfseries\normalsize\color{Navy}}
  {\thesubsubsection}{0.7em}{#1}[{\@addpunct{.}}]   
\makeatother
  
  \titlespacing*{\subsection}
  {0pt}{1.8ex}{0.8ex}

  \usepackage{titletoc}

\titlecontents{section}
  [0em]
  {\addvspace{0.3em}\sffamily\bfseries\color{Navy}}
  {\contentslabel{2em}}
  {}
  {\hfill\contentspage\vspace{-.3em}}

\titlecontents{subsection}
  [2.2em]
  {\vspace{-.3em}\small}
  {\contentslabel{2.5em}}
  {}
  {\titlerule*[0.5pc]{.}\contentspage}
  [\vspace{-.2em}]
  
  \makeatletter
  \titlecontents*{subsubsection}[2.2em]{\vspace{-.2em}\small\parindent=0pt\itshape\relax\S~}{\thecontentslabel\enspace}{}
	{{\@addpunct{,}} p.~\thecontentspage}[ \S ][.]
\makeatother

  \usepackage{setspace}
\usepackage{tikz}
  \usetikzlibrary{arrows.meta,decorations.text}
\usepackage{booktabs}
\usepackage{caption}
\usepackage{array}

\usepackage{fancyhdr}
  \usepackage{lastpage}
	\usepackage[hang,flushmargin]{footmisc}

\usepackage[most]{tcolorbox}
  \definecolor{PQback}{HTML}{F2F4F7}
\newtcolorbox{pullquote}{%
  enhanced, breakable, sharp corners, boxsep=0pt, boxrule=0pt,
  colback=PQback, colframe=Navy, leftrule=2pt,
  left=10pt, right=10pt, top=8pt, bottom=8pt,
  fontupper=\itshape\color{DarkGray}}

  \newcommand{\key}[1]{\emph{#1}}
  \newcommand{\keys}[1]{{\raggedleft{\ttfamily\footnotesize\color{LightGray}#1}\par}\vspace{.2em}}
  \newcommand{\ksep}{\,\textperiodcentered\,}

  \newcommand{\titleblock}[5]{%
    \begingroup
    \centering
    {\sffamily\bfseries\LARGE\color{Navy} #1\par}
    \vspace{0.45em}
    {\sffamily\large\color{Navy!72} #2\par}
    \vspace{0.75em}
    {\color{Navy}\rule{\linewidth}{0.9pt}}\par
    \vspace{0.75em}
    {\sffamily\normalsize\color{DarkGray} #3\par}
    \vspace{0.2em}
    {\sffamily\footnotesize\color{LightGray} #4\par}
    \vspace{0.45em}
    {\ttfamily\footnotesize\color{LightGray} #5\par}
    \vspace{1.4em}
    \endgroup}
  \newenvironment{abstractblock}
    {\begingroup\small\color{DarkGray}
     \setlength{\leftskip}{1.6em}\setlength{\rightskip}{1.6em}
     \noindent{\bfseries\color{Navy}Abstract.}\ \ignorespaces}
    {\par\endgroup\vspace{1.2em}}
\usepackage[
    backend=biber,
    style=alphabetic,
    sorting=nyt,
    maxbibnames=99,
    giveninits=true,
    doi=true,
    url=true,
    isbn=false,
    backref=true,        
    backrefstyle=three,  
]{biblatex}
		\AtBeginBibliography{\small\vspace{-.5em}}
		\defbibenvironment{bibliography}
		{\list
			{\bfseries [\printfield[labelnumberwidth]{labelalpha}]}
			{\setlength{\labelwidth}{0pt} 
				\setlength{\leftmargin}{0pt} 
				\setlength{\labelsep}{4pt}
				\setlength{\itemsep}{4pt}
				\setlength{\parsep}{2pt}} 
			}
		{\endlist}
		{\item}
		
\usepackage{xurl} 
  \definecolor{linkblue}{HTML}{1D4ED8}
  \definecolor{citegreen}{HTML}{047857}
\usepackage[%
  pdfauthor={Benjamin Collas},
  pdftitle={Fundamental Mathematics in the Age of AI -- The Residue, the Journey, and the Ecology},
  bookmarksnumbered=true,hidelinks=true,
  colorlinks=true,
    linkcolor=Navy,
    citecolor=citegreen,
    urlcolor=Navy,
]{hyperref}
  
  \usepackage{enotez}
		\setenotez{reset = false, counter-format = Alph,backref = true,list-style=paragraph, list-heading ={\section*{#1}} }
		
    \DeclareInstance{enotez-list}{mypara}{paragraph}{format = \small}
      \setenotez{list-style = mypara}

\begin{document}

\titleblock
  {Fundamental Mathematics in the Age of AI}
  {The Residue, the Journey, and the Ecology}
  {Benjamin \textsc{Collas}}
  {Research Institute for Mathematical Sciences, Kyoto University, Kyoto 606-8502, Japan}
  {bcollas@kurims.kyoto-u.ac.jp \ \textperiodcentered\ \today
  }

\vspace{2cm}

\begin{abstractblock}
  Large language models have begun refuting long-standing conjectures and, for a few thousand dollars of tokens, solving long-open problems (OpenAI, August 2026). The introspection this has prompted about the future of mathematical discovery is now well under way, and the anxiety accompanying it legitimate -- but both, we claim, are attached to the wrong loss. 
  
  What machines now produce is the countable part of mathematics -- theorems, proofs, refutations -- which was always the \emph{residue} of the work, and not its product. The distinction is old, and not originally economic: a result can be taken in its finished essence, or in the operations that engendered it. The product is human understanding: not a stock of results but a collective, hard-won way of deciphering the world and acting upon it\endnote{From this perspective, it should be no surprise that \emph{the book of nature ``is written in the language of mathematics''} (Galileo 1623) -- or closer to us, that the advanced theory of topos is applied to AI (2025). For a debate on the end game of fundamental mathematics or its utilitarianism, we refer to a 1830 comment from Carl Gustav Jacobi (result of Humboldt's mathematics as reine Wissenchaft as a source of freedom and personal realisation) to Joseph Fourier (the Polytechnicien of the French revolution, for which Mathematics as a state instrument is required) -- see \cite{JCor75}\begin{quotation}\small\itshape
    Fourier avait l'opinion que le but principal des mathématiques était l'utilité publique et l'explication des phénomènes naturels ; mais un philosophe comme lui aurait dû savoir que le but unique de la science, c'est l'honneur de l'esprit humain, et que sous ce titre, une question de nombres vaut autant qu'une question du système du monde. 
  \end{quotation} The maxim comes a few months before the opening of the Liverpool-Manchester steam railway (Sept. 1830), which set off the Railway Mania of the 1840s, and, with it, the globalisation of the industrial revolution.}. The two are arcs of a single loop: understanding tells us where to look; looking produces the residue; and taking it up again, one journey at a time, is what rebuilds the shared understanding. Machines are strong on the countable arc, absent from the one that feeds it. \emph{The peril is to leave the loop open.} AI did not create the confusion between residue and product; it has called a bluff long on the books, driving the cost of the residue towards zero and making the scarce thing visible at last.

  A new instrument makes a new way of working before it makes a new result. The pressing questions are therefore institutional: who can check an announced result, whoever announces it; what the work and the training become for the next generation of researchers; and whether the one thing that cannot be mass-produced -- the journey that nourishes a shared understanding -- continues to be funded. Mathematics, we argue, is uniquely placed among the sciences on the first -- a proof answers to no one's permission -- and uniquely exposed on the other two: teaching cannot go on as before, and no collective position yet exists; and the journey has never had a price our institutions knew how to pay. The decision is ours.
\end{abstractblock}

{\small  MSC (2020): {\ttfamily 01A80 (Primary) 68T01, 68V20, 00A35 (Secondary)}}

\newpage
  \thispagestyle{empty}

  \subsubsection*{Forewords}
    {\small \itshape These notes set down a reflection of some three years' standing. Their writing was occasioned by a conversation with Parmy \textsc{Olson} (Bloomberg) at the end of July 2026, and occupied the 27th and the 31st; they can be read as a companion paper to \cite{OLS26}. The OpenAI announcement of 1 August arrived as the first version was being finished, and, to my surprise, fitted the argument without altering it. They draw on personal experience as a researcher, and on exchanges with mathematical communities, private companies, and representatives of governing agencies. 
    
    Reference to economic principles should not be read as a claim of expertise, but as an attempt to reach across to those who have it. Political economy and access are not the frame of this essay, though it meets them where its argument does: openness and verification (\S~\ref{sec:asymmetry}), who can appropriate what and who provides it (\S~\ref{sub:extRiv}), dependency and sovereignty (\S~\ref{sec:structures}), entry cost (Coda). The choice is deliberate -- each would demand another register, and a competence this text does not claim.

    One question is declined outright: what machine intelligence is. This essay does not ask whether a machine understands; it asks whether the loop is closed, and by whom -- understanding being, as Part~3 argues, not a property anything possesses but a passage someone has to make again (\S~\ref{sec:stake}).

    \emph{Why a working mathematician should write this at all, rather than get back to work?} The seventh recommendation of the Leiden Declaration makes it a duty: mathematicians ``have a responsibility to support serious science journalism and to engage in public discourse to explain and contextualize artificial intelligence-assisted methods and results'' \cite{leiden2026}. The questions addressed in this note will be settled by whoever turns up to discuss them.
    }
\newpage

  \setcounter{tocdepth}{3}
  \tableofcontents

  \hfill

\subsubsection*{Acknowledgements}
{\small  The author thanks \href{https://www.bloomberg.com/authors/AVYbUyZve-8/parmy-olson}{Parmy \textsc{Olson}}, whose questions occasioned these notes; \href{https://cornebise.com/julien/}{Julien \textsc{Cornebise}}, whose good-humoured and productive disagreement with an earlier version sharpened the arguments that follow; \href{https://ochigame.org/}{Rodrigo \textsc{Ochigame}} for warm exchanges on AI, Lean, and understanding in May 2026; \href{https://www.kurims.kyoto-u.ac.jp/~motizuki/}{\textsc{Mochizuki} Shinichi} for regular exchanges on how the Math-AI frontier could be shaped; \href{https://pro.univ-lille.fr/pierre-debes/}{Pierre \textsc{Dèbes}} for comments on a previous version of these notes. He also thanks \href{https://patrickshafto.com/}{Patrick \textsc{Shafto}} (expMath) for checking the passages describing his work, and is grateful to \href{https://webusers.imj-prg.fr/~ariane.mezard/}{Ariane \textsc{Mézard}}, \href{https://www.math.univ-paris13.fr/~wittenberg/}{Olivier \textsc{Wittenberg}} and \href{https://research-portal.uu.nl/en/persons/johan-commelin/}{Johan \textsc{Commelin}} for comments on Version $1$ of the manuscript, and the latter for bringing \cite{OST10} to his attention. Thanks also to colleagues in several countries, within the \href{https://ahgt.math.cnrs.fr/}{``Arithmetic and Homotopic Galois Theory'' project} and beyond, for many conversations on the future of scientific discovery.
}%


\newpage
\mbox{}\medskip
{\raggedleft
\begin{minipage}{.6\textwidth}
  \small\itshape
  [\ldots] le résultat n'est plus considéré dans son essence si pure soit-elle
  devenue, mais dans la simple position d'un résultat qui résulte, par suite
  indication des modes qui l'ont engendré.\par
  \vspace{0.35em}
  \mbox{}\hfill\upshape\footnotesize
  -- J.~\textsc{Cavaillès}, \emph{Sur la logique et la théorie de la science}
  \cite{CAVA47}\footnotemark
\end{minipage}
\footnotetext{``[\ldots] the result is no longer considered in its
  essence, however pure that may have become, but in the simple position of a
  result that results -- an indication, thereby, of the modes that engendered
  it.'' We follow the reading of \cite{CN01}, for whom Cavaillès thematises,
  behind the dialectic of concepts, a \emph{mathematical experience}: the
  mathematician's experience, in practice, of the internal dynamics of the
  concepts themselves -- neither a psychology of discovery nor a
  socio-historical determination of mathematical content.}
\par}

\vspace{1.4em}

\section{The economy of the residue}\label{part:residue}

{\itshape Three claims, in order: that a conservative field is moving faster than it has in generations, and on tools it does not own; that what machines now produce is a real but narrow part of mathematics, and was always the residue of the work rather than its product; and that the economics of the field had already mistaken the one for the other, long before any model was trained.}

\subsection{A conservative field, moving fast}\label{sec:conservative}

\keys{Thurston\ksep human understanding\ksep Leiden Declaration\ksep ICM~2026\ksep means of production}

Academia -- and fundamental mathematics in particular -- was built as a place for producing sound knowledge in the service of human understanding\footnote{\textbf{Premise.} Throughout, ``understanding'' means human understanding, and "the field" the human community that reads, checks and teaches. Whether a machine could hold understanding of its own, or constitute an alternative ecosystem with its own criteria, is a real question and not this essay's.\label{FT:one}}: not a stock of results, but a way of deciphering the world and of acting upon it. This is the position on Mathematics of Thurston (Fields Medal 1982): he opens by refusing the obvious question -- not ``how do mathematicians prove theorems?'', but \emph{``how do mathematicians advance human understanding of mathematics?''} (\S~1) -- and closes where this essay starts:

\begin{quotation}\itshape
  [The proposed division of roles into ``speculation'' and ``proving''] only perpetuates the myth that our progress is measured in units of standard theorems deduced. This is a bit like the fallacy of the person who makes a printout of the first 10,000 primes. \key{What we are producing is human understanding.} We have many different ways to understand and many different processes that contribute to our understanding.\hfill\textup{-- W. Thurston \cite[\S~5]{thurston1994}}
\end{quotation}

The intangibility of the product, or rather its uncountability, makes us, researchers in mathematics, conservative in adopting new tools and new ways of working -- no instrument can quickly prove that it helped. It is therefore worth registering how fast the field is currently moving.

The most valuable tool of a researcher in fundamental mathematics is their voice, so the explosion of discussion is itself the measure. Over the past two years researchers  have engaged strongly with AI tools:
\begin{enumerate}
  \item \emph{In practice:} writing, investigating hypotheses, producing examples, paving the field -- and working groups forming inside academia for the purpose\footnote{For instance the HANA project at RIMS, Kyoto University, on AI-assisted formalisation for anabelian geometry and related topics -- members include \textsc{Hoshi} Yuichiro, \textsc{Mochizuki} Shinichi and the author -- and its partner group LANA.}.
  \item \emph{In print:} the \emph{Bulletin of the AMS} gave AI and mathematics two special issues in 2024 \cite{AMSAI24}, and individual reflections have followed -- Akshay \textsc{Venkatesh} (Fields Medal 2018) -- see~\cite{VEN24};
  \item \emph{In community:} the Leiden Declaration \cite{leiden2026} was drawn up in under a year, and five mathematicians have set out a joint position in the \emph{Notices} \cite{CJOTV26};
  \item \emph{In public:} the ICM~2026 programme includes a public lecture by Terence Tao (Fields Medal 2006), \emph{Mathematics in the Age of AI}, see \cite{icm2026};
  \item \emph{In funding:} establishment of dedicated funding structures such as DARPA's \emph{expMath} \cite{expmath}, which explicitly aims at accelerating progress in pure mathematics\footnote{expMath oversees multiple projects to shape the interface between AI and a broad range of mathematical disciplines -- from number theory to partial differential equations and homotopy theory -- combining fundamental research with societal impact.} -- and with private capital, in new structures that hire researchers and publishing at pace \cite{AXI26}.
\end{enumerate}
A field that changes on the scale of two generations has reorganised its conversation in a semester.

    Absorbing new modes of producing knowledge is, in one sense, ordinary work for researchers in mathematics: we did it with computers in the 1990s (\S~\ref{subsub:Comp}), and we are doing it with mass collaboration and formalisation since \emph{circa} 2020. \key{What is not ordinary this time is that the means of production belong to someone else.}

\subsection{What is being measured}\label{sec:measure}

\keys{Counterexample search\ksep choice of definitions\ksep ten problems for \$2{,}000\ksep Goodhart's law\ksep what a benchmark measures}

Mathematics is often presented as the ideal measure of AI progress: a proof is right or wrong, where prose is a matter of taste. The measurement is real, but it is worth asking what it measures.

\subsubsection{Counterexamples, and the constructive half}

Machines have so far produced results of two kinds: counterexamples to standing conjectures, and -- more recently, and beyond refutation -- solutions to open problems (\S~\ref{sec:buys}). We take them in turn.

The first kind is counterexample search -- killing conjectures, most recently, and among many, \emph{the Jacobian conjecture}, by Levent Alpöge on July 20, 2026 \cite{RHB26}. This is real, and we do not minimise it: counterexample search is the region of mathematics most amenable to machine methods. It is also the least representative of a mathematician's work, and it already illustrates, structurally, the missing interface between AI companies and academia (\S~\ref{sec:work}).

\begin{pullquote}
Disproving a conjecture is like pruning a branch that grew in the wrong direction. It does not tell us in which direction to look; it is not an end result. A counterexample, like any result, is worth the horizon it opens.
\end{pullquote}

The other half of the work -- the constructive half -- is deciding which definitions are worth making and which objects deserve names. That half is what makes a field rather than a pile of results: \key{a theorem is only as good as the notions it is stated in.} As of today, no system yet performs that choosing\footnote{And the field is moving fast. Beyond this question, what machine's developing interaction with fundamental mathematics do create is a new frontier to investigate: the shape of the line defining the area where they perform.}. Pruning accelerates a field only if someone is still choosing which branches were worth growing\footnote{Interestingly, it is precisely how the conversation has evolved once the counter-example to the Jacobian conjecture has been obtained: ``what does the counter-example tell us on the statement?'' We isolated this example for that reason -- a turning point in the mathematics and AI dynamic.}.

\subsubsection{What \$2{,}000 really buys}\label{sec:buys}
  The second kind goes further. OpenAI has recently reported solutions to what they present as ten long-open problems, obtained for roughly \$2{,}000 of tokens \cite{openai2026}\footnote{Ten successes at roughly \$2{,}000 apiece -- but out of how many attempts? What about the \$200--500 million development cost for each version of a frontier model with a six-month life expectancy \cite{SCH26}? The announcement does not say: the figure prices a success, a \emph{marginal cost}, not the search that found it.}, and -- the part worth registering -- published them in good order: manuscripts prepared with human help, arguments shipped with a Lean certificate, and, following a certain interpretation of the Leiden Declaration \cite{leiden2026}, an explicit refusal to claim human authorship. The same announcement concedes the limit of operating apart from academia: these questions ``cannot be answered by a technology company alone.''

  Good order, however, is not yet good practice\footnote{OpenAI's solution to \emph{3.~Non-sofic groups} (existence) is a case in point: the narration rests on an unrefereed preprint (cited as ``[Kun19]'') containing typographic errors, and on what experts regard as a strong hypothesis. The existence of a non-sofic group was, moreover, widely-accepted (see \cite{FFF26}).\label{ft:OAI3}}. The manuscript -- a narrative of the discovery -- names no one: no co-author to write to, none of the integrators and writers behind the ``human help''. The letter of the standard is met while its point -- a human interface on which the field can build -- is left empty; we return to this \emph{open loop} in \S~\ref{sec:work}.

  It is worth being exact about what the price tag measures. 
  
  \begin{pullquote}
    What \$2{,}000 buys is not a theorem; it is evidence of capability at the tasks the benchmark selected -- and that evidence holds only for as long as the community that certifies it still functions.
  \end{pullquote}
  Even that evidence is softer than it looks. Much of what a corpus contains has never been noticed -- the connections standing to be realised in a literature nobody can read, what Bessis calls the ``overhang'' \cite{BES26} -- and a large overhang certainly weakens the claim behind ten long-open problems solved: how many were still, in the relevant sense, open? The non-sofic group is a case in point (Footnote~\ref{ft:OAI3}).

  Note that the claim is about \emph{independence}, not about the indispensability of researcher in mathematics -- see \S~\ref{sub:extRiv}. What is mainly at issue is not the correctness of the proof -- which is or isn't independently of the reader -- but the inference from them to a capability.

   This is Goodhart's law\footnote{``When a measure becomes a target, it ceases to be a good measure'' \cite{goodhart1975}. Mathematics has watched this happen before, from the inside: publication counts and impact factors became targets, and ceased to measure what they once tracked. The novelty is not the mechanism but the scale, and that the party optimising the measure now sells it.}, and here it is not a metaphor: the measure is itself the product being sold.

   Generative systems are especially exposed to this issue: they are ungrounded, and the party with the financial incentive to optimise the measure is the party producing it.

\subsubsection{Capability rests on context}
  On capability itself, our assessment is plain. A year ago, LLMs behaved in mathematics like badly educated students. They have since been very well trained, and can now identify wrong directions of thought and good directions for problem solving. But their contribution stays narrow in scope\footnote{``[C]ompared to human papers, LLM ideas are disproportionately concentrated around bridge-like opportunities and synthesis methods [...] with a narrower range'' \cite{CZC26} -- a visible consequence of the mechanism that consists in harvesting the overhang.} and still rests on context previously built -- the definitions, the libraries, the accumulated understanding. Whether that limit is conceptual or merely a matter of engineering, we do not know, and we would not bet against scale: nothing we know of prevents a system from proposing its own definitions -- see Footnote~\ref{FT:one}. 

  Working mathematicians are already wondering where the reach can stops:

  \begin{quotation}\itshape
      There are many who will argue that AI is going to come up with genuinely new ideas that people have never thought of before. I don't know that I believe that [\ldots] But the ideas are somehow already there.\\[0.2em]
      \mbox{}\upshape\hfill -- K.~Ono \cite{ONO26}
  \end{quotation}

  The reason is structural. A model is placed to digest a corpus of ideas and results no person can read -- a researcher covers under 0.1\% of it in a working lifetime -- and to harvest what is already standing in it; at that scale, harvesting and creating become indistinguishable to the non-expert, for whom the research frontier is invisible in any case\footnote{The traffic runs the other way too. For Klowden and Tao, mathematics is a ``sandbox'' in which the impact of AI can be observed \cite[\S~3]{KT26}; for us, it is more than that. \emph{Confronted with a practice whose object is to make ways of understanding explicit, a machine shows the shape of whatever understanding it does have} -- which parts are within its reach, and where the reach stops. Whether that amounts to an understanding of its own is not this essay's question (Footnote~\ref{FT:one}); but it is in mathematics that such a question would become visible at all. AI came out of mathematics; mathematics is what can take its measure.}. And harvesting a corpus is not producing one.\emph{An age of gleaning rather than of discovery?}

   Suppose the doubt is misplaced, and that machines do come to choose the notions. The concession costs the argument nothing. If machines come to choose the notions, \key{what changes is who chooses -- not that choosing is the scarce act.} Someone must still recognise a good choice as good. One might reach for a trained notion of ``goodness'', but a trained ranker optimises an existing proxy, while mathematics favours a notion on which one can build -- and that is settled downstream, not at ranking time. It is that recognition, that \emph{judgment}, not the generating, that this essay is about.

    Whether or not they one day generate their own context, the question will remain that of the human--machine interface: its nature, and its quality (\S~\ref{part:journey}).

  \subsection{The economics were already wrong}\label{sec:economics}

    \keys{Residue and product\ksep public good\ksep externality\ksep the unfunded ecology\ksep mathematical capitalism}

    \subsubsection{No one ever worked out how to buy it}

      The economic story about mathematics was wrong before any of this. We are funded as if theorems were the product, and we sell teaching as the justification. Both are residue\footnote{For a cultural explanation of the overvaluation of the theorem -- Hardy's curse: an honour code that priced theorems and discounted everything else -- we refer to \cite{BES26}: the same starting point in Thurston and a similar diagnosis, but a different cause and a different remedy. For Bessis, the curse must be undone by saying publicly what mathematics is; for us, it is a failure of the economic and institutional contract, which leads to the questions of who verifies, who trains, and who pays. A mechanism closer to ours, because it is economic rather than cultural, is the currency metaphor Klowden and Tao attribute to Minhyong \textsc{Kim}: mathematicians must ``accumulate some credibility `currency', by proving difficult new mathematical results, before they can `afford' to `spend' that currency on speculative activities, such as formulating conjectures or philosophizing about the broader consequences of a result'' \cite[\S~4.6]{KT26}.}. The product is human understanding -- and no one ever worked out how to buy it\footnote{The temptation here is the one Avigad names: ``Threatened by such an argument, mathematicians may retreat to aesthetics: many of us do mathematics not because it is useful, but because we enjoy it, just as we enjoy literature and art'' \cite{AVI25}. We decline the retreat. The claim of these pages is not that mathematics is a pleasure society should indulge, but that it is the one such practice able to verify the claims of a trillion-dollar industry without that industry's permission (\S~\ref{sec:asymmetry}): supporting it is self-interested before it is cultural. Classical music cannot audit anyone.}. 

      The distinction is not that the theorem is worthless -- it is the instrument of transmission -- but that what it transmits is \emph{a capacity to think that must be rebuilt in each reader}; the writing, the seminar, the refereeing are that rebuilding, and none of it is what the theorem states. And the residue does not only transmit: it orients. A conjecture is pure residue -- countable, publishable, and not even proved -- yet its whole office is to make a horizon visible, so that someone knows where to look next; that is the arc Figure~\ref{fig:loop} draws, from the residue through understanding to the choice of the next question.

      One case is worth having in full, because the institution counts it only as overhead, and because the researchers themselves rarely notice it either: \emph{the preparation of a funding application.} Written as science rather than as a bid, it yields a research programme ready to be carried out -- by the funder's own survey, half the programmes of unsuccessful applications are carried out all the same, on other money. \emph{The application was refused; what it produced was kept.}

      The figures are not small -- the mathematical sciences are credited with \pounds495 billion of gross value added in the United Kingdom, and with about one job in eight there and in France alike\endnote{Orders of economic magnitude -- with the caveat that such studies measure the economic activity mathematics \emph{enables}, not the value it produces. In the United Kingdom, the mathematical sciences were credited with \pounds208 billion of gross value added and 2.8 million jobs for 2010 -- \emph{about 16\% of GVA and 10\% of employment} \cite{deloitte2012} -- and with \pounds495 billion and 4.2 million jobs, 13\% of UK employment, for 2023 \cite{acadmathsci2024}. In France, the corresponding study puts the share of \emph{GDP impacted by mathematics at 18\% and the jobs at 3.3 million, 13\% of salaried employment}, on 2019 data, up from about 16\% and 2.4 million on 2012 data \cite{CNRS22}. The trend is the interesting part; so is the fact that \emph{none of these studies attempts to price what this essay calls the product}.} -- and yet the membrane between what we understand and the economy has always been close to impermeable.

      \key{AI has not created this problem. It has called the bluff}, because it can mass-produce the residue (Table~\ref{tab:residue}). What it cannot produce is the reason anyone wanted the residue in the first place.

\begin{table}[t]
\begin{threeparttable}
\caption{\emph{Residue and product.} The claim of this essay is the last row: AI drives the marginal cost of the left-hand column towards zero, and thereby makes the right-hand column visible, for the first time, as the scarce thing.}
\label{tab:residue}
\small
\renewcommand{\arraystretch}{1.25}
\begin{tabular}{@{}>{\raggedright\arraybackslash\ttfamily}p{3cm}%
                  >{\raggedright\arraybackslash}p{5cm}%
                  >{\raggedright\arraybackslash}p{4.9cm}@{}}
\toprule
 & \textbf{The residue} & \textbf{The product} \\
\midrule
What it is
  & Theorems, proofs, counterexamples; courses and credentials\tnote{a}
  & Understanding: the choice of notions, judgment, the horizon\tnote{a} \\
Countable
  & Yes
  & No \tnote{b}\\
Reproducible
  & Always was -- by any trained mathematician
  & Never -- remade, not copied\\
Mass-producible by machine
  & Yes, and increasingly so
  & No \\
What pays for it
  & Grants, publication counts, teaching
  & Nothing, directly \\
What AI changes
  & \itshape Its marginal cost falls towards zero
  & \itshape Its scarcity becomes visible \\
\bottomrule
\end{tabular}
\begin{tablenotes}
        \item[a] A journey is what deposits both columns at once: a countable proof and the horizon it opens.
        \item[b] Counting the people who have understood a notion counts instances of transmission, not the understanding itself: the collective grasp of a notion is not the sum of individual graspings.
    \end{tablenotes}
\end{threeparttable}
\end{table}

      \subsubsection{Calling the bluff}\label{sec:bluff}
          The bluff was, in fact, already on public record. In 2019, Japan's ministries of economy and education published a joint report under the remarkable title \emph{The Coming Era of Mathematical Capitalism} \cite{meti2019}, which names the country's top three priorities for science:
          \begin{quotation}\itshape
            ``mathematics, mathematics, and mathematics.''\\[-.2em]
              \mbox{}\upshape\hfill -- METI and MEXT (2019)
          \end{quotation}
          It observes that the corporate race in AI ``requires human resources with a higher knowledge of mathematics,'' and asks -- explicitly, and without answering -- what ``management or social system'' would suit such an era. The question of our social contract is thus not a private anxiety of mathematicians; it is an open item on the books of governments. 
          
          Five years on, the same ministry is funding an attempt at an answer: a CREST area opened in 2024 under Motoko \textsc{Kotani}, \emph{Creation of Mathematical Foundation for Prediction and Control}, whose stated aim is to ``build a new social infrastructure'' on the basis of mathematical prediction \cite{JST24} -- a structure rather than a project (\S~\ref{subsub:UnitStruc}).

          Behind both the question and the attempt is a single change: the membrane between what we understand and the economy, always close to impermeable, has become permeable in both directions -- mathematics is drawn into the economy as a priced input, and economic demand now reaches back into what mathematicians are funded to do. \emph{That is, the traffic is priced.}

          From a broader societal perspective -- with the economy as one of its parts -- what crosses the membrane is the portion of the good that carries a price; what still stays behind is the rest -- a public able to read and integrate what its machines produce, the slow formation of people, a cultural shift. When does the greater part begin to cross?

\subsubsection{The externality, and the non-rivalry of the corpus}\label{sub:extRiv}

  The new arrangement gives the old imbalance a precise name.

\begin{pullquote}
  Responsibility for correctness is the cheap half, and it is now machine-dischargeable; the expensive half -- reading, situating, repairing, teaching, turning a correct proof into something that can be built upon -- is left to us, researchers, and takes years.
\end{pullquote}

That is an \key{externality} in the ordinary economic sense: extremely well-capitalised consumers; an unfunded ecology of researchers who keep understanding alive (``maintainers'', says the software idiom -- a word that fits the libraries and undersells both the people and the collective intelligence), still needed, by OpenAI's own admission (\S~\ref{sec:buys}) -- see also the AI training requirements of \S~\ref{sub:extRiv}; and infrastructure that rots quietly. 

But an externality is only an argument if it names an instrument; ours would be structural rather than fiscal, and we return to it in \S~\ref{sec:structures}.

Two things in particular are consumed and neither is paid for. \emph{And here the economics is, unusually, on the argument's side:} what is known of the returns to basic research places them in people and their accumulated understanding, not in the results they publish -- the residue and the product, under an older name; non-rivalry is why the corpus, not the theorem count\footnote{Here we depart from \cite[\S~5]{KT26}, who read the question distributively -- who gains, who is displaced -- and rehabilitate the Luddites as rational actors defending a livelihood against a technology that took it. The distinction drawn here cuts elsewhere. The \emph{residue can be appropriated:} priced, sold, mass-produced, held. \emph{The product cannot:} Understanding is non-rival, and no corporation can own a share of it -- which is why the danger to it is not capture but defunding, and why the remedy is institutional rather than distributive.}, is the asset. \emph{The first is the corpus:} not the literature, but the digested, exposited, taught mathematics the models are trained on, which took centuries to build and which does not deplete by being read. But \key{non-rivalry holds of the stock, not of its provision}: what depletes is the pipeline that keeps producing and digesting it -- the finite attention of the people who expose, referee, canonise and teach. It is a different kind of good, and the one actually under pressure\endnote{Economists have treated basic research as a public good since Nelson \cite{nelson1959} and Arrow \cite{arrow1962}: non-rival, non-excludable, and therefore privately under-supplied -- the standard argument for funding it publicly. \emph{Non-rivalry is also the engine of endogenous growth theory}, in which the accumulation of ideas rather than of capital sustains growth \cite{romer1990} -- which puts the corpus, and not the theorem count, at the centre of the economic story. The review literature adds something more useful here. Surveying the evidence, Salter and Martin \cite{salter2001} find that the benefits do not travel mainly through the published result: they sort them into \emph{six channels, of which \textbf{(1)}~codified knowledge is only one}, the others being \textbf{(2)}~trained people, \textbf{(3)}~new instruments and methods, \textbf{(4)}~networks, \textbf{(5)}~problem-solving capacity, and \textbf{(6)}~new firms. 

That is, in an economist's vocabulary, the distinction this essay draws between the residue and the product, the limitation of problem-solving results, and the requirement for funding new integrated structures. \emph{A non-specialist's question follows:} if most of the return travels through people rather than through results, what becomes of the return when the results are free and the people are not funded?

Elinor \textsc{Ostrom} (Nobel Prize in Economic Sciences, 2009) answers by re-typing the good, correcting Samuelson's dichotomy: ``rivalry of consumption'' becomes \emph{subtractability of use}, and a fourth type appears -- the \emph{common-pool resource}, subtractable like a private good and hard to exclude from like a public one \cite{OST10}. The corpus is not that good; its provision is. What depletes is the pipeline, drawn on a finite pool of expert attention from which no one can be excluded -- see also Footnote~\ref{ft:ostrom} of \S~\ref{sec:structures} for some structural consequences.}.

\emph{The distinction is not academic, and the companies discovered it empirically:} Practicioners report that AI training on the published literature and on forums like MathOverflow underperformed -- though both are already digested for human readers -- and AI companies found they had to commission new corpora, written by mathematicians for the purpose. For LLMs as well, the digested layer is not a by-product of publication; it has to be made, by people, deliberately.

\emph{The second is the measure} -- what \$2{,}000 buys is evidence of capability, and that evidence is only worth something while the community certifying it still functions (\S~\ref{sec:measure}). This is not a guild claiming to be indispensable: \emph{it is the ordinary requirement of any science that a measure be certified independently of the party it measures.} Kernels check proofs -- imperfectly, and only proofs: whether a statement says what was meant (\S~\ref{sec:asymmetry}), and whether ten instances warrant a claim about a system, no machine certifies.

Which raises the question to put to the ecosystem rather than to the models: can it produce understanding, economically and culturally, with no researcher in the loop, when its targets are whatever can be measured? 

{\itshape An oracle answers; a science explains -- one should not trade the second for the first without noticing.}

\section{The audit and the structures}\label{part:institutions}

{\itshape From what is produced to what can be required. Mathematics can verify without asking anyone's permission, which is rarer than it sounds, and it can therefore put questions to whoever announces a result. What it cannot do alone is build, fund, and hold open the structures on which that capacity depends.}

\subsection{The asymmetry, and the audit}\label{sec:asymmetry}

\keys{Verification without permission\ksep formalisation\ksep Lean and Mathlib\ksep training data\ksep attribution}

    \subsubsection{Verification without the producer's permission}

      What follows from the means of production belonging to someone else (\S~\ref{sec:conservative}) is not only a question of which interface gets built, but of scale: \emph{an asymmetry of manpower and of compute that no university can match.} This is serious for a field with a culture of openness and the cheapest science there is: a person, some paper, access to a library\footnote{\emph{And some coffee!} Following the polymath Paul Erdös' apocryphal ``A mathematician is a machine for turning coffee into theorems.''}... some mathematics!

      Mathematics, however, is the one field where that asymmetry has a technical answer rather than a negotiated one. Every other discipline that wants to check an AI company's claims must ask the company for access; the weights and the compute are theirs, and access can be withdrawn quietly, without anyone breaching an agreement. We do not have to ask (Table~\ref{tab:asymmetry}).

      \begin{pullquote}
        A proof -- written or formal -- can, and must, be checked independently of the prover: it has always been the rule of the field. No API key, no compute parity, no permission.
      \end{pullquote}

      The volume of machine production is a genuine difficulty, and this is where formalisation enters: machine-checked proofs, open source, carried by communities such as those behind Lean and its library Mathlib \cite{ML20}, or Rocq (formerly Coq), Isabelle and Agda. It is the one form of refereeing that scales with machine output, and it makes hallucination detectable -- at the level of the proof. \emph{But refereeing was never only error-detection, nor the certification of some absolute truth:} it is also how a field digests a result -- how one person's argument is read, situated, and made part of what everyone can build on. A kernel scales the first office and not the second\footnote{The same limit can also be reached from the side of the text: a proof carries a ``smell'' -- the impression an experienced reader forms before checking a single step -- and, around its deductive core, a ``penumbra of heuristic, empirical, or metamathematical reasoning'' that tells one why the argument works and whether it extends -- see \cite{KT26}. Trained on correctness alone, a machine can thus produce an ``odorless'' proof. Flat and mute, we would say: verified, and with no relief a reader can steer by, telling us nothing about why it works or how far it reaches. This ``smell'' is our \emph{judgment} (\S~\ref{sec:work}), turned toward reading rather than toward the choosing of definitions. What they find missing from the text is what this essay counts as the labour nobody funds.}.

      At the level of the statement, matters are different: no kernel can check that a formal statement says what the mathematics meant. The accident report of \cite[\S~5]{LS25} is instructive -- the authors relate discovering, mid-project, that a kernel-accepted statement was not the theorem they had meant to prove; a Lean evaluation convention, under which functions silently return ``junk'' values outside their domain of definition had changed the content of the statement, while leaving every proof valid -- see also \cite[\S~4.4]{KT26}. 

      But formalisation can also work the other way: in economic theory, formalising Aumann's agreement theorem -- fifty years old in 2026 -- made precise what it means for two parties to hold ``the same priors'', a hypothesis whose subtleties had gone unexamined \cite{ONO26}.

      Either way, the same conclusion: whether a statement says what was meant is not something a kernel settles, and formalisation can at best expose it. That is a question for people -- which is what the next section asks.

    \subsubsection{Our asset is their raw material}

        There is an irony here worth stating plainly. The property that makes formalised mathematics auditable -- that a proof can be checked by a kernel, with no human in the loop -- is exactly the property that makes it valuable as training material for LLMs. 

        The Leiden Declaration says so without euphemism: what makes mathematics attractive to general-purpose AI development is that ``the correctness of formalized proofs can be checked automatically, without the need for human oversight,'' which yields ``an effectively unlimited source of feedback for training'' \cite{leiden2026}. \key{The same property that lets us audit them lets them mine us.} 
        
        This is not an argument against formalisation. It is an argument for noticing that our defensive asset is also somebody's raw material, and for asking on what terms the library is being read.

    \subsubsection{Why do we cite?} 
        Attribution is a different matter, and \emph{deeper than an acknowledgement of debt.} Why do we cite? Not to pay persons. 
        
        We cite prior works as lighthouses: they mark where understanding advanced -- the progress, the techniques, the context that cast light and what it lit -- over a landscape still being charted. Models are trained on the published commons and, as the same Declaration observes, ``frequently return outputs that do not properly cite the human works they synthesize'' -- on data often obtained by ``systematically exploiting licenses and access arrangements that were not made with artificial intelligence in mind'' \cite{leiden2026}. No kernel adjudicates attribution, and here we are exactly as exposed as every other discipline. 

        What is lost is not credit: \key{an output without lighthouses flattens the landscape}, and a landscape without relief offers no journey -- nothing to steer by, nothing worth travelling to. Three of the questions below establish whether a claim is sound; only the fourth asks whether the field can build on it.

\begin{table}[t]
\caption{\emph{Who can check whom.} Mathematics is the only row whose verification does not require the producer's cooperation, and the only one with a form of refereeing that scales with machine output. The ``yes'' holds at the level of the proof; at the level of the statement it does not.}
\label{tab:asymmetry}
\footnotesize
\setlength{\tabcolsep}{4pt}
\renewcommand{\arraystretch}{1.25}
\begin{tabular}{@{}>{\raggedright\arraybackslash\ttfamily}p{2.0cm}%
                  >{\raggedright\arraybackslash}p{2.6cm}%
                  >{\raggedright\arraybackslash}p{3.0cm}%
                  >{\raggedright\arraybackslash}p{2.6cm}%
                  >{\raggedright\arraybackslash}p{3.0cm}@{}}
\toprule
 & \textbf{A claim rests on} & \textbf{Who holds the means of verification} & \textbf{Checkable without the producer's permission?} & \textbf{What scales with machine output} \\
\midrule
Fundamental mathematics
  & A proof
  & Anyone holding the statement and the argument
  & \textbf{Yes}
  & Machine-checked formalisation (Lean, Mathlib) \\
Empirical sciences
  & Data and an experiment
  & Whoever holds the apparatus, the sample, the cohort
  & In principle -- at the cost of replication
  & Nothing comparable: replication does not scale \\
Machine learning itself
  & A benchmark score
  & The company: weights, data, compute
  & No -- access is granted, and can be withdrawn quietly
  & Nothing: evaluation is itself model-mediated \\
\bottomrule
\end{tabular}
\end{table}

\pagebreak
  \subsection{Four questions for an announced result}\label{sec:questions}

    \keys{Disclosure\ksep who chose the problem\ksep independent evaluation\ksep First Proof\ksep EMS Code of Practice}

     \subsubsection{Whoever makes the announcement} 
      Hence the questions we should put to any announcement of a mathematical result, whoever makes it -- AI company or individual -- written or formal:

      \begin{enumerate}\itshape 
        \item Is the complete argument public -- the understanding, not only the residue?
        \item Who chose the problem, and the terms in which it is stated -- and when?
        \item What independent checking process has been applied?
        \item Who closes the loop -- is there a named interlocutor, and terms under which the field can read, question, attribute, and build on the result?
      \end{enumerate}

      The four bind an individual exactly as they bind a laboratory. The fourth, however, answers itself in one case and not in the other -- someone announcing their own result is the interlocutor, and can be written to; an announcement issued by a company -- as in the case of the solutions to ten long-open problems \S~\ref{sec:buys} -- frequently has no one behind it at all, which is the open loop of~\S~\ref{sec:work}.

      In the formal case\footnote{It may be of interest to note that mathlib Lean library is not only a technical instrument but a culture with its own norms: the reference library's own paper presents ``the social organisation that has led us here'' \cite{ML20}.} this also becomes precise: \key{who wrote the formal statement, and was it fixed before the model started?} If one system produced both statement and proof, the verification establishes mostly \emph{internal consistency} -- ground for suspicion rather than reassurance about the claim itself.

    \subsubsection{The third question, institutionalised}
      Something answering to the third question already exists. First Proof \cite{firstproof} is an academic-led nonprofit that benchmarks AI systems on research-level mathematics, its editorial and governing boards drawn from Stanford, Berkeley, Harvard, Yale, Columbia, Chicago, UT~Austin, and EPFL/Imperial College London.
      
      In its ``First Batch'' project \cite{FPRF}, it evaluated problems in a range of fields, under controlled conditions, with expert refereeing for correctness \emph{and} for exposition: of the ten problems it set, seven were solved at publication-level quality, at a compute cost between \$10 and \$1{,}000 apiece. The shape of the announcement echoes \S~\ref{sec:buys}; what differs is who ran it.

      The score is not the interesting output. The evaluation was designed and run by the people qualified to referee it, rather than reported by the party being evaluated -- \key{the third question, institutionalised.}

    \subsubsection{The standard already exists}
      None of this is a demand invented for the occasion. The European Mathematical Society's Code of Practice \cite{ems2012} already instructs mathematicians not to make public claims of new theorems unless full details can be provided in a timely manner. 
      
      And the extension of that standard beyond the profession is no longer a private view: the Leiden Declaration calls on collaborations between mathematicians and industry to abide ``at minimum, by the standards we expect of our colleagues'' \cite{leiden2026}. \key{We are asking no one -- company or individual -- for more than the standard mathematicians bind themselves to.}

      Questions and standards, however, only bind where the two sides actually meet (the fourth question). What remains is to build, and to fund, the places where they do -- the structures, to which we now turn.

\pagebreak      
  \subsection{The structures are the story}\label{sec:structures}

    \keys{Interfaces\ksep dependency at scale\ksep open code and publication\ksep evaluation and corpus\ksep sovereignty\ksep fund structures, not projects}

      \subsubsection{The two sides are not symmetric}

        Interfaces are now being built from two sides, and they will outlast any individual result: \textbf{(i)}~From the companies: OpenAI's \emph{ChatGPT for Academic Researchers} offers free model access to a hundred thousand scientists \cite{openai2026}, Google DeepMind's \emph{AI for Math Initiative} opens its most advanced tools to five partner institutions\footnote{Google Deepmind's partners: Imperial College London (UK), Institute for Advanced Study (US), Institut des Hautes Études Scientifiques (FR), Simons Institute for the Theory of Computing (US), Tata Institute of Fundamental Research (IN) -- see \cite{GDP25}.} -- access at scale on one side, access by selection on the other; \textbf{(ii)}~From academia: First Proof, expMath, Mathlib, and a growing set of instruments built by the community rather than borrowed from outside\footnote{For example: \url{https://www.erdosproblems.com}; Optimization Constants in Mathematics; Competitions of the SAIR type, which ``explores the frontier of AI and science through open challenges, benchmarks, and community collaboration''.}.

        The two sides are not symmetric, and the asymmetry cuts both ways. \emph{Free access at scale is dependency at scale, not a research partnership:} access is granted, and can be withdrawn without anyone breaching anything (\S~\ref{sec:asymmetry}), only now with a hundred thousand people downstream of the decision. In the other direction, two things sit on our side of the table and not on theirs. 
        
        The first is the evaluation:

          \begin{pullquote}
          A capability claim is worth what its certification is worth, and the certification is ours -- designed, run and refereed by the community, see \S~\ref{sec:asymmetry} -- The asymmetry, and the audit.  
          \end{pullquote}
          
          The second is the corpus, which academia does not merely hold but produces and keeps legible -- the definitions, the libraries, the expositions, and the people who write them. \key{Each side is short of exactly what the other holds.} 

          That is a negotiating position, and a better one than the current mood suggests -- a mood that presses for close collaboration and calls it symbiosis. It holds only if we notice what is actually being settled now: the structures, not the individual results.

        \subsubsection{Interface-structures: the form is not the practice}\label{subsub:Anchor}
          Between these two sides, however, a \emph{third kind of structure is emerging: interface-structures}, neither corporate laboratory nor university department -- {\ttfamily Axiom Math} or {\ttfamily MathInc} among them, a \emph{company} in the older sense of the word\footnote{A group of people with a shared interest who mean to make something happen in the world.} rather than in the sense these pages have so far had to give it. They hire doctoral students, early-career and senior researchers; they build new interfaces of knowledge, integrating AI capabilities, for sound and fast progress in understanding. What they need from academia is recognition, funding and partnership -- an alternative to sending our people to the companies, not an outpost of them. \emph{That anchor is what makes them interface-structures rather than laboratories:} without it, the incentive is to capture the recognition and leave the corpus-work to the community that has to live with it.

          The form, however, is not the practice. A structure of this kind can produce work a field can build on, or it can produce volume the field must then clean up: formalisations no one can read (or rather digest), results with no way into the literature (or no way to contextualise). \emph{Which of the two it produces depends on the anchor} -- representatives of academia with the standing to say what would count as usable, and answerable for having said it.  
          
          It depends, too, on what leaves the structure. Open code, open publications, and open access to new tools for results that anyone can ingest and digest. For two out of three, that is already the practice of the instrument the community has built\footnote{This forestalls the divide \cite[\S5.3]{KT26} anticipate. It is, besides, ordinary scientific practice -- open code and open data are what make a result reproducible, and reproducibility is the same asset as the audit of \S\ref{sec:asymmetry}.} -- and it is \emph{what keeps an interface-structure an interface rather than a gate.}
        
          Accepting such structures as partners is not a courtesy; it is the condition under which what they produce can be built on at all.

      \subsubsection{The unit to fund is a new type of structure}\label{subsub:UnitStruc}
          If there is an instrument here, it is complementary to the ones we have and not a rival to them\footnote{\label{ft:ostrom}``Complexity is not the same as chaos.'' In support of a multi-core structure network that gathers the efforts of a multitude of autonomous contributors -- of which the \emph{Lean Focused Research Organization} \cite{leanfro} and the \emph{Mathlib Initiative} \cite{mathlibinit} are examples -- we refer to Ostrom's economic notion of \emph{polycentricity}: technical efficiency becomes higher when many autonomous direct producers coexist with a small number of shared indirect ones \cite{OST10}, see also \S~\ref{sub:extRiv}.}. That instruments of this kind can be made well is not in doubt: the European Research Council funds a programme and a person rather than a deliverable, and a generation of European mathematicians owes to it the freedom to choose its own questions. It did not price the product -- no one has -- but it stopped requiring the residue in exchange.

          What no agency produced is expMath, Axiom Math, MathInc or Mathlib, which emerged outside the academic frame -- a gap in the catalogue of instruments, not a fault in any of them. The agencies have since begun to answer in their own form, anchored in academia: the NSF's Institute for Computer-Aided Reasoning in Mathematics (ICARM) opened in 2025 \cite{ICARM}, and the CREST area of \S~\ref{sec:bluff}.
      
          Designing a structure has not, until now, been the mathematician's business: the social contract gives us teaching and research, and delegates structuring to agencies. This is changing from both ends: the agencies are building institutes, and mathematicians are now among the founders of the new interface-structures. The direction can therefore be stated from the inside, as it should be.
      
          Two interfaces have a claim on such instrument. The newest, where \emph{academia and AI companies genuinely meet rather than transact}, has still to be built: the new institutes are anchored on the academic side, and it is the crossing that is missing. That crossing is \emph{the natural place for the externality of \S~\ref{sec:economics} to be internalised} -- the ecology funded where the corpus and the evaluation are actually produced, the people, not only the libraries they maintain. The oldest is education, which produces the value of the next generation rather than of the next quarter -- including the capacity to use these tools well -- and which therefore needs money on a clock that is not the quarter's.
      
          Both are questions of sovereignty and of culture as much as of economics. A country that stops training the people who can read and integrate what the machines produce has outsourced rather more than its compute. 

      \key{The unit to fund is the structure}: projects and individuals are countable, which is precisely why they are the ones we do.

\section{The ecology of the journey}\label{part:journey}

  {\itshape What all of this does to the work itself, and to the people who do it: the loop that constitutes the profession, the training of those about to enter it, the mismatch between the speed at which results are produced and the time understanding takes to form, and the one thing in the field that was never reproducible. These are questions of morale as much as of production.}

\pagebreak
  \subsection{The loop is the profession}\label{sec:work}

    \keys{The discovery loop\ksep judgment\ksep human-machine interface\ksep open loop\ksep doctoral training}

    \subsubsection{The fear, and the two cuts}

      A widely shared fear -- vividly put in a recent essay by a young mathematician \cite{hampshire2026} -- is that the profession survives only as a taste-labelling layer: paid ``to optimize the weights for `interesting', `beautiful', anything you like,'' permitted to comment and appraise, forbidden to create. 

      The same suspicion arrives at the other end of a career. 
      \begin{quotation}\itshape
        If my success as a mathematician relied only on my ability to learn techniques that somehow could be put together to prove a theorem, then maybe that was actually automatable, and maybe I mistook all of that hard effort for something maybe it wasn't.\\[-0.2em]
      \mbox{}\upshape\hfill -- K.~Ono \cite{ONO26}
      \end{quotation}

      What he had been doing instead, he says, was
      \begin{quotation}\itshape
      I was learning how mathematics fits together to help me become a mathematician who can ask these questions and participate in the discovery.
      \end{quotation}

      In both, the same displacement: what the profession counted turns out to be the part that automates, and what it never counted -- learning how mathematics fits together -- \emph{is what made the mathematician}. And in both, the same unspoken premise: that appraising and discovering are separable jobs, one of which can be kept while the other goes.

\begin{figure}[t]
\centering
\begin{tikzpicture}[
    >={Stealth[length=2.2mm,width=1.6mm]},
    every node/.style={font=\small},
    lp/.style={->,semithick,black!80},
  ]
  \def\R{2.75}

  \node[align=center] (und) at ( 90:\R) {\bfseries Understanding};
  \node[align=center] (que) at (  0:\R) {Choice of question,\\choice of definitions};
  \node[align=center] (inv) at (270:\R) {Investigation};
  \node[align=center] (res) at (180:\R) {Residue:\\theorems, proofs,\\counterexamples};

  \draw[line width=5pt,black!12] (inv) to[bend left=22] (res);
  \node[align=center,font=\small\itshape,black!65] at (247:1.45)
       {what machines\\now do well};

  \draw[lp] (und) to[bend left=22] (que);
  \draw[lp] (que) to[bend left=22] (inv);
  \draw[lp] (inv) to[bend left=22] (res);
  \draw[lp] (res) to[bend left=22] (und);

  \path (que) to[bend left=22] node[pos=.5,inner sep=0pt] (cutin)  {} (inv);
  \path (res) to[bend left=22] node[pos=.5,inner sep=0pt] (cutout) {} (und);
  \foreach \c/\a in {cutin/315, cutout/135}{%
    \begin{scope}[shift={(\c)},rotate=\a]
      \fill[white] (-0.24,-0.16) rectangle (0.24,0.16);
      \draw[semithick] (-0.22,-0.13) -- (0.22,-0.03);
      \draw[semithick] (-0.22, 0.03) -- (0.22, 0.13);
    \end{scope}}
  \node[align=center,font=\footnotesize\itshape] at (310:4)
       {\emph{Input:}\\label and weight};
  \node[align=center,font=\footnotesize\itshape] at (135:4.45)
       {\emph{Output:}\\check and clean up};
\end{tikzpicture}
\caption{\emph{The loop is the profession.} Understanding tells us where to look; looking is what produces understanding. Machines are now strong over one quarter of the circumference -- and not over the quarter that feeds it. The taste-labelling layer is this same loop cut at the two points where the human meets the machine.}
\label{fig:loop}
\end{figure}
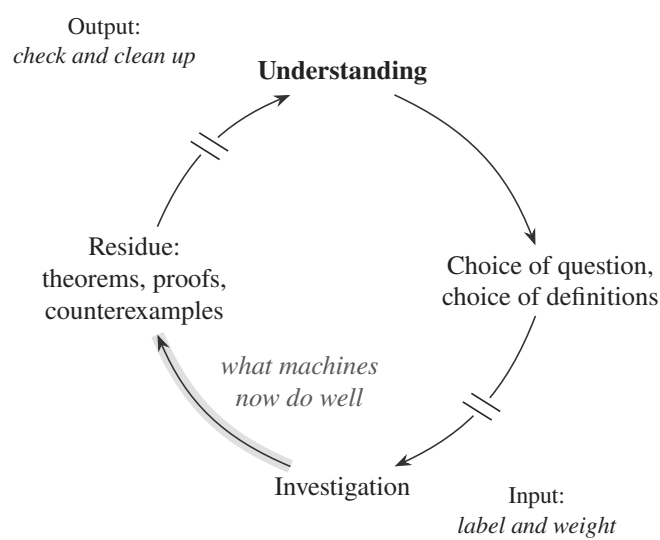

      The fear admits of a mechanical answer. Discovery and the weighting of ``interesting--beautiful'' are not two jobs but two moments of a single loop: understanding tells us where to look, and looking is what produces understanding. \key{The loop is the profession.} Cut it, and the damage arrives at both ends at once: at the \emph{output}, checking and cleaning up after the machine; at the \emph{input}, labelling and weighting with no discovery of one's own (Figure~\ref{fig:loop}). 

      The loop is drawn for one person; the collective understanding is not a larger loop but the ground the loops share: a notion belongs to the field when it can be used and re-walked by any of us. The residue is what travels between loops; the journey is what cannot.

\emph{How does one label \emph{where to look}, once detached from the understanding that alone tells us where to look?}

\subsubsection{The open loop}

The ten problems of \S~\ref{sec:measure} are a live specimen of that cut, and of what falls through it. At the input, there is no co-author. There is a published narration of the model's reasoning, which one can read but not question: nobody to write to, nobody to invite to a seminar, nobody who can be asked why this direction and not another. A narration is a record of a path, not a partner on it.

At the output, the contribution to understanding is real but in its weakest form: what Thurston calls a collection of ``answers'', where what people wanted was understanding \cite[\S~1]{thurston1994}. The result is correct, documented -- and still to be converted. That conversion, which is reading, situating, exposing and teaching, is years of work; it is ours; and nothing in the announcement funds it\footnote{In the formal case, this conversion is what Kontorovich calls ``canonization'': reworking a verified development into library mathematics general -- reusable and coherent enough that others can build on it, see \cite{KONT26}. It is the digested layer of \S~\ref{sub:extRiv} being made, deliberately, by people -- the word for the labour, not for what is produces.}. 

        So the proofs are correct and \key{the loop is left open} -- and an open loop is not a slower loop, it is a different object.

      \subsubsection{The loop kept whole: a new practice}

        The same loop, kept whole, points the other way. AI is \emph{a formidable trainer of intuition for reaching new research frontiers:} a doctoral student can test an intuition and build new understanding in afternoons where it used to cost months -- and so can any established researcher investigating outside their own field. Learning that one asked the wrong question was the most expensive lesson in a mathematical education; it is now nearly free.

        The change is not only one of tempo but of the nature of the work itself. \emph{Ingesting and digesting} -- reading and writing -- have long been dissociated acts, separated by months and often by years; the AI tool brings them into contact, so that one questions while reading and tests while writing. \emph{A new instrument makes a new way of working before it makes a new result.} The same instrument works at the other scale: a result can now reach, and be questioned by, more people -- and other people -- than any seminar could hold (Coda).

        What practice maintains, and a labelling layer cannot, is judgment -- whose deep form, in fundamental mathematics, is the choosing of definitions (\S~\ref{sec:measure}).

        \begin{pullquote}
          ``Knowing where to look'' is not a stock that can be harvested; it is maintained by practice. 
        \end{pullquote}

        Severed from practice, the layer degrades: the arrangement the fear imagines fails on its own terms.

        Until now, nobody funded question-asking; it produces nothing you can count. AI may change that: once the countable part costs nothing, the scarce part finally carries a price and acquires a value our socio-economic system had so far never assigned it. 

\subsection{Students get closer to the research frontier}\label{sec:StCloseRes}
  \keys{Doctoral training\ksep pre/post-training\ksep division of labour\ksep computer algebra, a precedent\ksep no collective position}

  {\itshape ``What is a mathematician trained to become, once the countable half of the work no longer needs one?''}

  The question is not hypothetical. It is being answered already, month by month, in the thesis topics proposed and the positions offered to those starting out -- and it is answered first for the people with the least say in it.

  \subsubsection{Not at the beginning of the loop}
    \key{The danger is not the tool but the routing.} The next generation is being pushed toward the two countable ends of the pipeline -- the ends that are easy to fund and easy to measure -- and those are precisely the stretches that train no judgment.

    The AI labs have names for the two ends: ``pre-training'' and ``post-training''\footnote{The terms are the industry's own: pre-training on a corpus, post-training by human correction of what the model produces. The mathematics student appears at both ends -- as a source in the first, as a corrector in the second, and as a mathematician in neither.} -- feeding the machine a corpus, then cleaning up after what it returns. Put a doctoral student at either end and one has built the taste-labelling layer of \S~\ref{sec:work} in advance, out of the people who were supposed to go round the loop in their turn. Such a layer does not degrade, as \S~\ref{sec:work} has it: \emph{it never forms.}

    On this point, we stand in opposition to the ``increasing division of labor'' of \cite{KT26}:
    \begin{quotation}\itshape
      One can envision an increasing division of labor in the future of mathematical research: [...] any given mathematician may increasingly specialize in just a few aspects of this process, for instance focusing on utilizing AI assistants to prove results as directed by some more senior member of a research group, or on using the most recent literature [...] to propose new directions of inquiry.\\[-0.2em]
      \mbox{}\upshape\hfill -- T.~Klowden and T.~Tao \cite[\S~4.6]{KT26}
    \end{quotation}
    Specialization at each side of the loop is the loop cut at both ends, and made into a career structure. 
    
    \begin{pullquote}
    A mathematical laboratory worth building is one in which each member is fully integrated and goes round the whole loop -- not one in which each is assigned an arc.
    \end{pullquote} 

    The objection is not to formalisation. It is a genuine competence, and the students who acquire it are not wasting their time: the objection is to formalisation as the whole of a \emph{formation} -- the French word for what \emph{training} has lost: preparing someone to digest the competences no one can yet name, rather than equipping them with the ones we can.

    A student who has never wandered the landscape has no way of learning that a question was the wrong one.

  \subsubsection{Classical training first, then contact, fast}

      Students still ask for a classical training, and are right to: definitions made by hand, proofs written out, the long slow reading. Reasoned use of the machine builds on that; it does not substitute for it.

    But the two ends of the pipeline are not the only place a student can be put. Contact must come early, supervised, and in three staged forms rather than one:
    \begin{enumerate}
      \item \textbf{Corpus.} Learning to read, and to ingest, at a scale no one could read before -- the machine as a way into a literature. This is where the overhang of \S~\ref{sec:buys} is met rather than described.
      \item \textbf{Acceleration of experience.} Testing an intuition in an afternoon instead of a term: learning that the question was wrong, cheaply.
      \item \textbf{Investigation.} The student's own question, with the machine inside the loop rather than at its ends.
    \end{enumerate}

    The order is not decorative. (2) without (1) produces fluency without reading; (3) without (2) produces confidence without calibration; (1) without (3) produces a reader who never asks. \emph{In each case, scientific progress is impoverished, if not entirely absent.}

  \subsubsection{The profession has already run this experiment: the 90's computer algebra}\label{subsub:Comp}

    Nothing here is a new anxiety: the last time a machine entered mathematics, the same two failures were available.

    When computer algebra arrived, the practice could have collapsed into \textbf{(a)}~a toolbox of computations and low-level programming, or \textbf{(b)}~a practice of pure first-order logic\footnote{The instruments arrived within a decade -- Maple (1982), PARI/GP (1983-1985), GAP (1986), Mathematica (1988), Magma (1993) -- and the journal \emph{Experimental Mathematics} was founded in 1992: the date at which the practice acquired a name, and a place to publish.}. Neither would have produced a mathematician. What it produced, when it worked, was \emph{a new paradigm of scientific thought} -- experimental mathematics -- and a generation that is able to think with machines rather than through them.

    The remedy is practical, not curricular. Judgment of that kind cannot be lectured into anyone: \emph{it is maintained by practice, and by nothing else.} Hence the arrangement to build: the student in front of the machine as early as possible, and supervised -- the one place where taste is called for, and therefore the one place it forms.

    Here we part from Bessis \cite{BES26}, for whom the frontier between teaching and research is sharper: in our view \emph{machine practice blurs it}, bringing students into contact with the mechanisms of research earlier than any curriculum could.

  Teaching as before is no longer possible, and no collective position exists -- no department, no learned society, no funder has one. The new practice is being developed, and valorized, in the new structures (\S~\ref{subsub:UnitStruc}); what is experimented there will percolate into academia -- which is how a practice becomes a \emph{formation}.

\subsection{Compression is not acceleration}\label{sec:rate}

\keys{Rate problem\ksep proof abundance\ksep company time vs.\ generation time\ksep Thurston on foliations}

\subsubsection{From proof scarcity to proof abundance}

AI compresses the stages of the work that can be counted -- from an open problem to an argument on the page. \emph{As the tools are use now}, they leave untouched the stage that follows: the one in which a result is read, situated, exposed, taught, and eventually becomes part of what the field simply knows -- the closing of the loop\footnote{Researchers are now doing exactly this for one of OpenAI's ten solutions: turning the narration into scientific progress (see \cite{MO513866} and Footnote~\ref{ft:OAI3}).}. That is how they are used, not what they can do: turned the other way, the same instruments lower the cost of entering a field instead of raising the volume it must absorb (Coda).

As Terence Tao (Fields Medal, 2006) states it plainly, if no ``suitable policy and cultural changes'' are made:

\begin{quotation}\itshape
  [W]e will transition from an era of proof scarcity to an era of proof abundance.\\[-0.2em]
\mbox{}\upshape\hfill -- T.~Tao \cite[slide 44]{icm2026}
\end{quotation}

Abundance is not, by itself, progress. The fact of a result travels cheaply; the horizon it opens does not. A thousand results nobody has travelled to will open little more than the fact of themselves.

\begin{pullquote}
  Compressing the countable stages does not accelerate understanding; it widens the distance between what exists and what is understood.
\end{pullquote}

\begin{figure}[h]
\begin{minipage}[t]{0.33\linewidth}\small\itshape
  \vspace{0pt}
  \caption{\emph{Time compression of generation is not, by itself, acceleration of understanding} (schematic). One turn of the loop, from an open problem to a result one can build on. The angle is each stage's share of the turn; the radius is how long the turn takes, so nearer the centre is faster; the thicker the ring, the more broadly the understanding is held.}
  \label{fig:compression}
  \medskip
  Machines compress the early stages in every case alike -- but only where the closing stage is carried out does the turn shorten. Left open, it never closes at all: compression of the early stages alone cannot reach understanding.
\end{minipage}\hfill
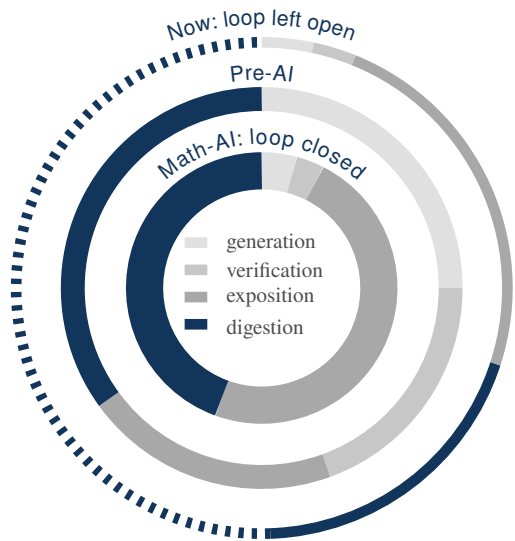
\begin{minipage}[t]{0.63\linewidth}
\vspace{0pt}
\centering
\begin{tikzpicture}[
    >={Stealth[length=2.0mm,width=1.5mm]},
    rC/.style={line width=14.0pt},   
    rB/.style={line width=9pt},   
    rO/.style={line width=4.0pt},   
    key/.style={font=\scriptsize,text=black!70,anchor=west},
  ]
  \def\Rc{1.55}   
  \def\Rb{2.50}   
  \def\Ro{3.25}   

  \draw[rC,black!12] ( 90:\Rc) arc ( 90:  75:\Rc);   
  \draw[rC,black!22] ( 75:\Rc) arc ( 75:  63:\Rc);   
  \draw[rC,black!34] ( 63:\Rc) arc ( 63: -110:\Rc);   
  \draw[rC,Navy]     (-110:\Rc) arc (-110:-270:\Rc);   

  \draw[rB,black!12] ( 90:\Rb) arc ( 90:   0:\Rb);   
  \draw[rB,black!22] (  0:\Rb) arc (  0: -70:\Rb);   
  \draw[rB,black!34] (-70:\Rb) arc (-70:-144:\Rb);   
  \draw[rB,Navy]    (-144:\Rb) arc (-144:-270:\Rb);  

  \draw[rO,black!12] ( 90:\Ro) arc ( 90:  78:\Ro);   
  \draw[rO,black!22] ( 78:\Ro) arc ( 78:  68:\Ro);   
  \draw[rO,black!34] ( 68:\Ro) arc ( 68: -18:\Ro);   
  \draw[rO,Navy]     (-18:\Ro) arc (-18: -88:\Ro);   
  \draw[rO,Navy,dash pattern=on 2pt off 3pt]
                     (-88:\Ro) arc (-88:-270:\Ro);   

  \path[decorate,decoration={text along path, text align={center}, raise=4.5pt,
        text={|\scriptsize\sffamily\color{Navy}|Math-AI: loop closed}}]
       (145:\Rc+0.20) arc (145:35:\Rc+0.20);
  \path[decorate,decoration={text along path, text align={center}, raise=2.5pt,
        text={|\scriptsize\sffamily\color{Navy}|Pre-AI}}]
       (145:\Rb+0.17) arc (145:35:\Rb+0.17);
  \path[decorate,decoration={text along path, text align={center}, raise=2.5pt,
        text={|\scriptsize\sffamily\color{Navy}|Now: loop left open}}]
       (145:\Ro+0.15) arc (145:35:\Ro+0.15);

  \begin{scope}[shift={(-1.02,0.60)}]
    \foreach \i/\c/\t in {0/black!12/generation, 1/black!22/verification,
                          2/black!34/exposition}
      {\draw[line width=4pt,\c] (0,-0.36*\i) -- (0.30,-0.36*\i);
       \node[key] at (0.42,-0.36*\i) {\t};}
    \draw[line width=4pt,Navy] (0,-1.08) -- (0.30,-1.08);
    \node[key,align=left] at (0.42,-1.15) {digestion};
  \end{scope}
\end{tikzpicture}
\end{minipage}
\end{figure}

    \subsubsection{What AI can shorten, and what it cannot}
      A field can move faster and understand less, at the same time. Economists have a name for the missing quantity: \emph{absorptive capacity}, the ability to use knowledge produced elsewhere, built only by doing research oneself \cite{cohen1990} -- which is why free access to results confers little on a community that has stopped producing the people who can read and integrate them (\S~\ref{sec:structures}). 

      Properly grounded, and held to the code of practice of the mathematical sciences (\S~\ref{sec:questions}), \emph{AI can however raise that capacity} -- the rate at which frontier work can be taken up -- by enhancing individual practice. The distinction matters: becoming standard is slow because it requires broad deliberative consensus, and AI does not speed deliberation. 
      
      What AI can shorten is \emph{each participant's climb to the point of being able to take part.}

\subsubsection{A mismatch of clocks}

  Underneath is a mismatch of clocks. \key{Company time is release-cycle time; understanding takes root on the time of a generation.} The two cannot be brought into phase by working harder at the fast end -- which is the only end anyone is currently paying to speed up.

  Nor is this a machine problem in origin: the precedent is at human scale and required no AI. Thurston's first subject was foliation theory, and his flood of results emptied it within a couple of years. Colleagues were advising one another not to go into foliations -- Thurston was ``cleaning it out'' -- and told him, not as a complaint, but as a ``compliment'', that he ``was killing the field''. He is explicit that the territory was not intellectually exhausted; what \emph{the rate of production destroyed were the \emph{ecological} conditions of the subject}, the ones under which anybody would want to keep thinking in it\endnote{Thurston's account, in \cite[\S~6]{thurston1994}, deserves quoting at length. Having ``fairly rapidly proved some dramatic theorems'' as a graduate student: \begin{quotation}\small\itshape An interesting phenomenon occurred. Within a couple of years, a dramatic evacuation of the field started to take place. I heard from a number of mathematicians that they were giving or receiving advice not to go into foliations -- they were saying that Thurston was cleaning it out. People told me (not as a complaint, but as a compliment) that I was killing the field. Graduate students stopped studying foliations, and fairly soon, I turned to other interests as well.\end{quotation} On the cause: ``I do not think that the evacuation occurred because the territory was intellectually exhausted -- there were (and still are) many interesting questions that remain and that are probably approachable.'' He names instead two ``ecological effects'': a high entry barrier, his results being ``documented in a conventional, formidable mathematician's style'' that ``depended heavily on readers who shared certain background and certain insights''; and the absence of anything left in it for anyone else -- ``More than the knowledge, people want personal understanding. And in our credit-driven system, they also want and need theorem-credits.'' 

  The remedy he adopted afterwards is the one this section is arguing for: ``In reaction to my experience with foliations [\ldots] I concentrated most of my attention on developing and presenting the infrastructure'', publishing definitions and ways of thinking while remaining ``slow in stating or in publishing proofs of all the `theorems' I knew how to prove'', so as to ``[leave] room for many other people to pick up credit.''}. 

  \key{A rate of production can destroy the conditions of understanding.} It is a rate problem before it is a machine problem -- and it now has a very much larger engine attached to it.

  \subsection{The decision is ours}\label{sec:stake}

    \keys{Priority vs.\ journey\ksep Library of Babel\ksep reproducibility\ksep authorship\ksep what cannot be mass-produced}

    \subsubsection{Priority, and the journey}

      A worry one hears from the next generation of researchers \cite{hampshire2026} has the shape of the Library of Babel -- Borges' library that contains every possible book, hence, somewhere on its shelves, every theorem, every proof, every refutation, already written. But that is already our condition: \key{every researcher assumes that any result they obtain could be obtained by another trained mathematician.}

      Is the experience of discovery being foreclosed? To answer, one must split ``discovery'' in two. One sense is priority -- first place in the world's ledger. That, machines will indeed make scarcer and noisier. The other is the journey -- the personal passage from confusion to understanding -- and it is indexed to a person, not to the ledger. 

    \subsubsection{Nothing is lost}

      The residue was always reproducible; \key{the journey never was.} Here the essayist quoted above puts it exactly right: ``human mathematics is Talmudic'' \cite{hampshire2026} -- the charge of learning comes from being conversant, across centuries, with the person who discovered. But then the conclusion inverts: if meaning travels through human authorship, the accumulated corpus -- the thing we spend our lives learning and teaching -- loses nothing to anything a machine proves afterwards. 

      Each of us adds a little to a great deal that was already there: almost everything one will ever understand was discovered by someone else, and that has never once made the understanding less ours. For understanding is collective only in the way it is kept alive: not stored but re-walked -- each journey rebuilds a piece of the shared corpus in another mind\footnote{A philosophically inclined reader will recognise a Hegelian movement here: mathematical thought reaching a reflexive stage at which the \emph{Geist} recognises itself in its own productions -- transcending any particular individual, medium or epoch, while subsisting nowhere but in those who take it up again.}. 

      That is why a machine's proof, which enters only the ledger, adds nothing to what the field actually understands: \key{the shared understanding lives only where someone has made the passage again.}

    \subsubsection{What is actually at stake}
      The channel of discovery is not being sealed. It has become crowded and noisy -- not a bad problem to have, and one we already build for: curated libraries, large collaborative interfaces, machine-assisted formalisation. 
      
      What is genuinely at stake is not on the ledger at all -- and it is a decision, not a fate.

      \begin{pullquote}
        The journey cannot be mass-produced -- and a collective understanding lives only where the journey is made. But it can be routed away, defunded, traded for what can be counted. That decision is ours, not the machine's.
      \end{pullquote}

      A decision of that kind shows itself in budget lines, in what a committee agrees to count, in the positions offered to people starting out -- and it is announced, plainly enough, in the collaborations struck with the new interface-structures (\S~\ref{sec:structures}).


      {\itshape The journey, then, is not a private consolation; it is the piece the rest of this essay defends: what closes the loop (\S~\ref{sec:work}), what the rate endangers (\S~\ref{sec:rate}), and what the audit and the structures of \S~\ref{part:institutions} exist to keep funded.}

\begin{figure}[h]
\centering
\begin{tikzpicture}[
    >={Stealth[length=2.2mm,width=1.6mm]},
    every node/.style={font=\small},
    lp/.style={->,semithick,black!80},
    lab/.style={font=\footnotesize\itshape,text=black!65},
  ]
  \def\R{2.7}

  \node[align=center] (mac) at ( 90:\R)
       {\bfseries Machines\\[0.1em]{\footnotesize produce}};
  \node[align=center] (ker) at (330:\R)
       {\bfseries The kernel\\[0.1em]{\footnotesize checks}};
  \node[align=center] (peo) at (210:\R)
       {\bfseries Researchers\\[0.1em]{\footnotesize integrate and choose}};

  \draw[lp] (mac) to[bend left=14]
      node[lab,pos=.5,right=3pt,align=left]
      {randomness} (ker);
  \draw[lp] (ker) to[bend left=14]
      node[lab,pos=.5,below=4pt,align=center]
      {rigidity} (peo);
  \draw[lp] (peo) to[bend left=14]
      node[lab,pos=.5,left=3pt,align=right]
      {judgment} (mac);
\end{tikzpicture}
\caption{\emph{The virtuous loop.} People choose which questions are worth the asking; machines produce at a rate no community could sustain; the kernel makes what they produce checkable without its author having to be trusted, which opens the door to whoever's argument checks. Researchers integrate what arrives into the understanding -- and only that arc, the return to judgment, closes the loop.}
\label{fig:virtuous}
\end{figure}
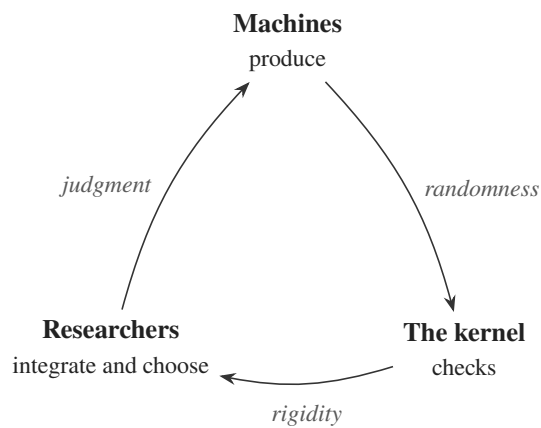

\section*{Coda}
\addcontentsline{toc}{section}{Coda}

  The models will keep producing mathematics, and mathematics will keep absorbing new modes of production, as it always has. The open questions do not belong to the machines: who checks, who trains, who pays -- and whether we continue to organise, and to fund, the one experience that was never reproducible in the first place.

  \emph{A virtuous loop nonetheless becomes available:} machines propose at a rate no community could sustain; formalisation makes the proposals checkable without their authors having to be trusted; and, the material now being machine-tractable, \emph{the same tools can be turned back on it} -- questioned, replayed, explained -- so that expertise becomes something a newcomer builds rather than a barrier already cleared (Figure~\ref{fig:virtuous}). And once the corpus once grounded by researchers, those tools let far more of them read, question and answer one another than any seminar could hold. 
  
  Read against Thurston's case, \emph{the virtuous loop lowers the entry cost of high-level knowledge} -- for a more inclusive and diverse practice.

  None of that is automatic: it holds only where the practice is built for it. Such tools can be turned outward or inward: they can prepare more people to take part in a conversation, or they can answer in its place. \emph{Which happens is settled by the structures, not by the model.} 
  
  The ``randomness'' from the LLMs, the ``rigidity'' from the kernel, and the judgment from people -- and the loop finally closes.

   \vfill

\begin{minipage}[t]{.45\linewidth}
  \vspace{0pt}
  \subsection*{Benjamin \textsc{Collas}}
  
  \medskip

  \small 
  Researcher at the Research Institute for Mathematical Sciences of Kyoto University (RIMS), working in arithmetic homotopy geometry -- a field that reconciles the rigidity of numbers with the flexibility of geometry. He leads the CNRS International Research Network \textsc{AHGT}, which links RIMS, the Universit\'e de Lille and the ENS-PSL, and brings together some $\sim$120 researchers in $\sim$15 countries.
  
  \medskip

  Email: \url{bcollas@kurims.kyoto-u.ac.jp}

  Web: \url{https://collas.perso.math.cnrs.fr}
\end{minipage}
\hfill
\begin{minipage}[t]{.45\linewidth}
  \vspace{0pt}
  \subsection*{Declaration on the use of AI}

  {\medskip \itshape
  The arguments, judgments and factual claims in this essay are the author's own. An AI assistant (Claude Fable~5 and Opus~4.8, Anthropic, accessed via Claude Code, July~2026) was used to draft and restructure English prose from the author's own writing and notes, for copy-editing, and editorial critique. 
  
  The author reviewed and revised all text, verified the sources, and takes full responsibility for the content.
  }  
\end{minipage}

\newpage 
{\addcontentsline{toc}{section}{Notes and References}
	\printendnotes
	}

  \subsubsection*{\sffamily References}

	\begin{multicols}{2}
		\printbibliography[heading=none] 
	\end{multicols}

\vfill

\subsubsection*{\sffamily Version history}
\begin{itemize}\small 
  \item 2026.08.29 -\emph{Publication.} ArXiv submitted (26 pages, 3 figures, 2 tables)
  \item 2026.08.28 - {\ttfamily Version 2.} New section on students and teaching; the loop kept whole as a new practice; the retreat to aesthetics declined \cite{AVI25}; the overhang of \cite{BES26}; the unit to fund as a new type of structure; the provision of the corpus as a common-pool resource, and polycentricity as the form the new structures take \cite{OST10}; declared scope on political economy and access; position with respect to \cite{KT26}; new sources: the 2024--2026 literature, K.~Ono (20 July 2026), Aumann on ``same priors''; ArXiv update (26 pages; 3 figures, 2 tables).
  \item 2026.08.13 - \emph{Publication.} Bloomberg published \cite{OLS26}; ArXiv submitted (17 pages, 2 figures, 2 tables);
  \item 2026.08.05 - \emph{Minor Additions.} The investigation on non-sofic groups, Figure~\ref{fig:compression};
  \item 2026.08.01 - {\ttfamily Version 1.}
\end{itemize}
\end{document}